\documentclass[11pt,a4paper]{article}
\usepackage{amsmath}
\usepackage{amsmath,amsthm}
\usepackage{xspace,amssymb}
\usepackage[colorlinks=true,allcolors=blue]{hyperref}

\theoremstyle{plain}
\newtheorem{theorem}{Theorem}

\newtheorem*{theorem*}{Theorem}
\newtheorem*{lemma*}{Lemma}
\newtheorem*{proposition*}{Proposition}
\newtheorem*{conjecture*}{Conjecture}
\newtheorem{fact*}{Fact}

\theoremstyle{definition}

\newtheorem{problem}[theorem]{Problem}

\newtheorem{remark}[theorem]{Remark}

\newtheorem*{definition*}{Definition}
\newtheorem*{question*}{Question}
\newtheorem*{example*}{Example}
\newtheorem*{remark*}{Remark}
\newtheorem*{remarks*}{Remarks}
\newtheorem*{exercise*}{Exercise}
\newtheorem*{assumption*}{Assumption}

\newcommand{\E}{\mathbb{E}}

\newcommand{\lesm}{\le_{\rm{sm}}}
\newcommand{\gesm}{\ge_{\rm{sm}}}
\newcommand{\leism}{\le_{\rm{ism}}}

\newcommand{\lec}{\le_{\rm{c}}}

\newcommand{\R}{\mathbb{R}}

\newcommand{\citeBB}{\cite{Blaszczyszyn_Yogeshwaran_2009}\xspace}

\newcommand{\HP}{\cite{Hu_Pan_2000}\xspace}

\newcommand{\cMSSS}{\cite{Muller_Stoyan_2002,Shaked_Shanthikumar_2007}\xspace}

\newcommand{\lei}{\le_{\rm i}}
\newcommand{\lecsm}{\le_{\rm csm}}
\newcommand{\leicsm}{\le_{\rm icsm}}

\begin{document}

\title{Ross's second conjecture and supermodular stochastic ordering\thanks{This preprint has not undergone peer review or any post-submission improvements or corrections. The Version of Record of this article is published in Queueing Systems, and is available online at https://doi.org/10.1007/s11134-022-09824-0}}

\author{Lasse Leskel\"a\footnote{
Department of Mathematics and Systems Analysis,
Aalto University,
Otakaari 1, 02015 Espoo, Finland.
}}

\date{7 February 2022}

\maketitle


\section{Introduction}


%


In 1978, about eight years before the first issue of \emph{Queueing Systems: Theory and Applications}, Sheldon M.\ Ross presented in \cite{Ross_1978} two famous conjectures about the classical single server queue subject to a randomly time-varying rate arrival rate $\lambda(t)$.
%
Intuition suggests that more autocorrelated arrival rate processes lead to more bursty arrival patterns, and therefore to larger workloads and longer waiting times. This question is naturally approached by modelling the arrival times using a Cox process (doubly stochastic Poisson process) in which the instantaneous arrival rate $\lambda(t)$ is a stationary and ergodic random process. 

\emph{Ross's first conjecture} in \cite{Ross_1978} states that the mean stationary waiting time in the system becomes smaller if the Cox arrival process is replaced by a homogeneous Poisson process with constant rate $\E \lambda(0)$. 
In 1981, Tomasz Rolski \cite{Rolski_1981} published an elegant proof this conjecture in full generality.  After five more years, Rolski \cite{Rolski_1986} also managed to prove a complementary upper bound, showing that the mean stationary waiting time in the system becomes larger if the intensity of the Cox process is replaced by a random constant $\lambda(0)$. These bounds display two extreme dependence structures:
\begin{itemize}
\item The lower bound corresponds to a fully averaged random environment and leads to completely independent (Lévy) dependence structure.
\item The upper bound corresponds to a frozen or quenched random environment, and leads to a maximally autocorrelated dependence structure.
\end{itemize}



\emph{Ross's second conjecture} in \cite{Ross_1978} concerns an arrival process with randomly time-varying intensity $\lambda_c(t) = \lambda(ct)$, parameterised by a modulation rate $c \in (0,\infty)$ and a stationary and ergodic baseline intensity $\lambda(t)$.  The parameter $c$ corresponds to a modulation of rate of a random environment, so that $c \to 0$ corresponds to freezing the environment to a (random) constant, and $c \to \infty$ corresponds to observing the environment at its ergodic average value.  Intuitively, a smaller modulation rate should correspond to a more autocorrelated intensity process, more bursty arrival point patters, and longer waiting times. Therefore, Ross conjectured that the mean stationary waiting time in the queue should be a decreasing function of the modulation rate $c$.


\section{Problem statement}

It appears that, after more than 40 years, Ross's second conjecture still remains not fully solved.
This motivates highlighting the following problem.
Consider a single server queue with with a Cox arrival process and IID service times, where the arrival rate $\lambda_c(t) = \lambda(ct)$ is defined in terms of a modulation rate $c \in (0,\infty)$ and a stationary and ergodic random baseline intensity $\lambda(t)$ for which the queueing system is stable.  Denote by $w(c)$ the mean stationary workload in the system.

\begin{problem}
\label{que:Ross2}
 Is it always true that
\begin{equation}
 \label{eq:Ross2}
 c_1 \le c_2 \quad \implies \quad w(c_1) \ge w(c_2),
\end{equation}
and if not, what is a simple necessary and sufficient condition for $\lambda(t)$ that is needed for \eqref{eq:Ross2} to be valid?
\end{problem}



\section{Discussion}

During past decades, the study of the above question has inspired the development of novel methods, e.g.\ stochastic dependence orders \cite{Muller_Stoyan_2002,Shaked_Shanthikumar_2007},  for analysing spatial and temporal autocorrelations in depth far greater than with pairwise correlation coefficients.  Key milestones partially answering the above question include:
\begin{itemize}
\item In the original article \cite{Ross_1978}, Ross proved \eqref{eq:Ross2} for a special case in which the baseline intensity is a two-state CTMC with one state being zero.
\item In 1991, Chang, Chao, and Pinedo \cite{Chang_Chao_Pinedo_1991} proved \eqref{eq:Ross2} in case where the baseline intensity is a finite-state CTMC with transition rates of form $q_{ij} = \alpha_i > 0$ for all $i \ne j$.
%
\item In 1998, Bäuerle and Rolski~\cite{Bauerle_Rolski_1998} established \eqref{eq:Ross2} in case where
the baseline intensity is a finite-state CTMC for which the transition rate matrix $Q$
and its time-reversal $Q^*$ with entries $Q^*_{ij} = \frac{\pi_j}{\pi_i} Q_{ji}$ are stochastically monotone with respect to the strong stochastic order on the real line.
\end{itemize}

As a consequence of latter milestone, we know that \eqref{eq:Ross2} holds whenever the baseline intensity $\lambda(t)$ is
a reversible and stochastically monotone CTMC, such as a two-state CTMC, or a birth-and-death process.
However, monotonicity in the strong stochastic order is a property of a dynamic system preserving the ordering of states \cite{Leskela_2010,Leskela_Vihola_2013,Massey_1987} and as such is not directly related to autocorrelation properties. If stochastic monotonicity precisely characterises CTMCs satisfying \eqref{eq:Ross2}, this would be a fundamental, and perhaps surprising, discovery.


Bäuerle and Rolski \cite{Bauerle_Rolski_1998} also showed that a sufficient condition for \eqref{eq:Ross2} is that $c \mapsto \lambda_c = (\lambda_c(t))$ is decreasing in the supermodular stochastic ordering, in the sense that
\begin{equation}
\label{eq:DecrSM}
 c_1 \le c_2 \implies \lambda_{c_1} \gesm \lambda_{c_2},
\end{equation}
where we recall that for stochastic processes $X \lesm Y$ means that the inequality
$\E \phi(X_{t_1}, \dots, X_{t_n}) \le \E \phi(Y_{t_1}, \dots Y_{t_n})$ holds for all integers $n \ge 1$, all time instants $t_k$, and all supermodular functions $\phi: \R^n \to \R$ for which the expectations exist \cite{Muller_Stoyan_2002,Shaked_Shanthikumar_2007}. Furthermore, any doubly stochastically monotone CTMC studied in \cite{Bauerle_Rolski_1998} satisfies \eqref{eq:DecrSM}. These discoveries have inspired impressive analytical works analysing supermodular regularity of Markov chains \cite{Hu_Pan_2000,Meester_Shanthikumar_1993,Miyoshi_Rolski_2004} and more recently also spatial random measures and point patterns \cite{Blaszczyszyn_Yogeshwaran_2009}. Nevertheless, the following problems appears are still unanswered in full generality to date.

\begin{problem}
\label{pro:NotMon}
Does there exist a stationary and ergodic CTMC $\lambda = (\lambda(t))$ which is not stochastically monotone but satisfies \eqref{eq:DecrSM}? 
\end{problem}

\begin{problem}
\label{pro:DecrSM}
Describe a simple necessary and sufficient condition for a stationary and ergodic CTMC $\lambda = (\lambda(t))$ under which \eqref{eq:DecrSM} holds.
\end{problem}

Problems \ref{pro:NotMon} and \ref{pro:DecrSM} are appealing because, outside the domain of stochastically monotone CTMCs, not much seems to be known about them even for state spaces of size 3. Furthermore, the above problems may also be formulated for discrete-time Markov chains, either by uniformising the CTMC, or by analysing chains generated by dampened transition probability matrices $(1-c)I + cP$ with $c \in (0,1]$ and $P$ representing a baseline Markov chain.  To the best of my knowledge, all of the above problems remain equally open also in the discrete-time setting.


\bibliographystyle{abbrv}
\bibliography{lslReferences}

\newcommand{\SortNoop}[1]{}\def\cprime{$'$}
\begin{thebibliography}{10}

\bibitem{Bauerle_Rolski_1998}
N.~B{\"a}uerle and T.~Rolski.
\newblock A monotonicity result for the workload in {M}arkov-modulated queues.
\newblock {\em J. Appl. Probab.}, 35(3):741--747, 1998.

\bibitem{Blaszczyszyn_Yogeshwaran_2009}
B.~B{\l}aszczyszyn and D.~Yogeshwaran.
\newblock Directionally convex ordering of random measures, shot noise fields,
  and some applications to wireless communications.
\newblock {\em Adv. Appl. Probab.}, 41(3):623--646, 2009.

\bibitem{Chang_Chao_Pinedo_1991}
C.-S. Chang, X.~Chao, and M.~Pinedo.
\newblock Monotonicity results for queues with doubly stochastic {P}oisson
  arrivals: {R}oss's conjecture.
\newblock {\em Adv. Appl. Probab.}, 23(1):210--228, 1991.

\bibitem{Hu_Pan_2000}
T.~Hu and X.~Pan.
\newblock Comparisons of dependence for stationary {M}arkov processes.
\newblock {\em Probability in the Engineering and Informational Sciences},
  14(3):299--315, 2000.

\bibitem{Leskela_2010}
L.~Leskel{\"a}.
\newblock Stochastic relations of random variables and processes.
\newblock {\em J. Theor. Probab.}, 23(2):523--546, 2010.

\bibitem{Leskela_Vihola_2013}
L.~Leskel{\"a} and M.~Vihola.
\newblock Stochastic order characterization of uniform integrability and
  tightness.
\newblock {\em Statist. Probab. Lett.}, 83(1):382--389, 2013.

\bibitem{Massey_1987}
W.~A. Massey.
\newblock Stochastic orderings for {Markov} processes on partially ordered
  spaces.
\newblock {\em Math. Oper. Res.}, 12(2):350--367, 1987.

\bibitem{Meester_Shanthikumar_1993}
L.~E. Meester and J.~G. Shanthikumar.
\newblock Regularity of stochastic processes: {A} theory based on directional
  convexity.
\newblock {\em Probability in the Engineering and Informational Sciences},
  7(3):343--360, 1993.

\bibitem{Miyoshi_Rolski_2004}
N.~Miyoshi and T.~Rolski.
\newblock Ross-type conjectures on monotonicity of queues.
\newblock {\em Aust. N. Z. J. Stat.}, 46(1):121--131, 2004.
\newblock Festschrift in honour of Daryl Daley.

\bibitem{Muller_Stoyan_2002}
A.~M{\"u}ller and D.~Stoyan.
\newblock {\em Comparison Methods for Stochastic Models and Risks}.
\newblock Wiley, 2002.

\bibitem{Rolski_1981}
T.~Rolski.
\newblock Queues with nonstationary input stream: {R}oss's conjecture.
\newblock {\em Adv. Appl. Probab.}, 13(3):603--618, 1981.

\bibitem{Rolski_1986}
T.~Rolski.
\newblock Upper bounds for single server queues with doubly stochastic
  {Poisson} arrivals.
\newblock {\em Mathematics of Operations Research}, 11:442--450, 1986.

\bibitem{Ross_1978}
S.~M. Ross.
\newblock Average delay in queues with non-stationary {P}oisson arrivals.
\newblock {\em J. Appl. Probab.}, 15(3):602--609, 1978.

\bibitem{Shaked_Shanthikumar_2007}
M.~Shaked and J.~G. Shanthikumar.
\newblock {\em Stochastic Orders}.
\newblock Springer, 2007.

\end{thebibliography}
\end{document}